\documentclass[12pt]{article}
\usepackage{amssymb,amsmath}
\usepackage[dvips]{graphics}

\numberwithin{equation}{section}

\newcommand{\Bialg}{\rm Bialg}

\newcommand{\Gr}{G r}

\newcommand{\Cbu}{C^{\bullet}}

\newcommand{\cCb}{\check{C}^{\bul}}
\newcommand{\cpa}{\check{\partial}}

\newcommand{\Der}{{\rm Der}\,}
\newcommand{\Hom}{{\rm Hom}\,}
\newcommand{\Alg}{{\rm Alg}\,}
\newcommand{\Coalg}{{\rm Coalg}\,}

\newcommand{\sgn}{{\rm s g n}}

\newcommand{\lan}{\langle}
\newcommand{\ran}{\rangle}

\newcommand{\lie}{{\bf lie}}
\newcommand{\comm}{{\bf comm}}
\newcommand{\ass}{{\bf assoc}}
\newcommand{\colie}{{\bf colie}}
\newcommand{\cocomm}{{\bf cocomm}}
\newcommand{\coass}{{\bf coassoc}}
\newcommand{\Ger}{{\bf e_2}}

\newcommand{\Bra}{{\bf Braces}}

\newcommand{\Hocomm}{Ho({\bf comm})}

\newcommand{\Hoger}{Ho({\bf e_2})}

\newcommand{\bul}{{\bullet}}

\newcommand{\al}{{\alpha}}
\newcommand{\la}{{\lambda}}
\newcommand{\io}{{\iota}}

\newcommand{\Om}{{\Omega}}
\newcommand{\Si}{{\Sigma}}
\newcommand{\si}{{\sigma}}

\newcommand{\ve}{{\varepsilon}}

\newcommand{\mgl}{{\mathfrak{gl}}}

\newcommand{\pa}{{\partial}}

\newcommand{\cB}{{\cal B}}
\newcommand{\cM}{{\cal M}}

\newcommand{\cC}{{\cal C}}

\newcommand{\cV}{{\cal V}}
\newcommand{\cO}{{\cal O}}

\newcommand{\bbR}{{\mathbb R}}

\newcommand{\bbK}{{\mathbb K}}
\newcommand{\bbH}{{\mathbb H}}

\newcommand{\bbF}{{\mathbb F}}

\newcommand{\La}{{\Lambda}}

\newcommand{\D}{{\Delta}}

\newcommand{\erarrow}{\stackrel{\sim}{\rightarrow}}

\newcommand{\brarrow}{\succ\rightarrow}

\date{}
\newtheorem{defi}{Definition}
\newtheorem{pred}{Proposition}

\newtheorem{teo}{Theorem}
\newtheorem{cor}{Corollary}

\title{The homotopy Gerstenhaber algebra
of Hochschild cochains of a regular algebra
is formal}

\author{
Vasiliy Dolgushev,
Dmitry Tamarkin, and
Boris Tsygan}

\begin{document}

\maketitle

\begin{abstract}
The solution of Deligne's conjecture on
Hochschild cochains and the formality of
the operad of little disks
provide us with a
natural homotopy Gerstenhaber
algebra structure on the Hochschild cochains of an
associative algebra.
In this paper we construct a natural chain
of quasi-isomorphisms of homotopy Gerstenhaber
algebras between the Hochschild cochain
complex $\Cbu(A)$ of a regular commutative algebra
$A$ over a field $\bbK$ of characteristic zero
and the Gerstenhaber algebra of multiderivations
of $A$\,. Unlike the original approach of
the second author based on the computation
of obstructions our method allows us to
avoid the bulky Gelfand-Fuchs trick and
prove the formality of the
homotopy Gerstenhaber algebra structure on
the sheaf of polydifferential operators on a smooth
algebraic variety, a complex manifold,
and a smooth real manifold.
 \end{abstract}

\section{Introduction}
Various proofs of Deligne's Hochschild cohomological
conjecture \cite{KS}, \cite{M-Smith}, \cite{Dima1}
\cite{Sasha} by
J.E. McClure, J.H. Smith, M. Kontsevich, Y. Soibelman,
A. Voronov and the second author show that the
complex $\Cbu(A)$ of Hochschild cochains of
an associative algebra $A$ is equipped with
a remarkable structure of a homotopy Gerstenhaber
algebra. The homotopy Lie part of this structure
coincides with the DGLA structure provided by
the Gerstenhaber bracket \cite{Ger} and the
corresponding commutative (up to homotopy) product
is cohomologous to the cup product.

If $A$ is the algebra of polynomials
$\bbK[x_1, \dots, x_n]$ over the field
$\bbK$ of characteristic zero
then it is not hard
to prove formality of the homotopy Gerstenhaber
algebra $\Cbu(A)$ by explicitly computing
the obstructions \cite{Dima11}, \cite{Dima1}.

In paper \cite{Gilles}
G. Halbout showed that the structure maps
of a $G_{\infty}$ quasi-isomorphism
for the Hochschild cochains of the algebra
$\bbR[[x_1, \dots x_n]]$ can be made relative
to Lie subalgebra $\mgl_n\subset \Der(\bbR[[x_1, \dots x_n]])$\,.
This result suggests that the Gelfand-Fuchs trick
\cite{CEFT},
\cite{Fedosov}, \cite{GF},  \cite{G-Kazh}, \cite{VdB},
\cite{Ye} could be applied to globalization
of the $G_{\infty}$ quasi-isomorphism.

In this article we construct a natural chain
of quasi-isomorphisms of homotopy Gerstenhaber
algebras between the Hochschild cochain
complex of a regular commutative algebra $A$ and
the Gerstenhaber algebra of multiderivations
of $A$\,. An important advantage of this
approach in comparison to the original
proof \cite{Dima11} and \cite{Dima1} of the second
author is that we propose an explicit construction
of a chain of quasi-isomorphisms while in
\cite{Dima11} and \cite{Dima1} only the
existence of the desired chain is proved.
This approach allows us to
avoid the bulky Gelfand-Fuchs trick
and prove the formality of the
homotopy Gerstenhaber algebra
structure on the sheaf of
polydifferential operators on a smooth
algebraic variety, a complex manifold,
and a smooth real manifold.

The organization of this paper is as follows.
In the introductory section we explain notation
and recall basic facts about (co)algebras and
(co)operads. In the second section we discuss
properties of the homotopy Gerstenhaber algebra
structure
on the Hochschild cochain complex induced by the
action \cite{M-Smith} of the singular chains
of the operad of little squares. In section $3$
we prove the main result of this paper,
theorem \ref{domik}, and in section $4$ we
give obvious generalizations of this theorem.
In the concluding section we discuss open
questions and further possible applications
of our results.

{\bf Acknowledgment.}
We are thankful to G. Halbout and M. Markl
for discussions.
D.T. and B.T. are supported by 
NSF grants. D.T. is also supported by the
research fellowship of A.P. Sloan.
The work of V.D. is partially supported the grant CRDF
RM1-2545-MO-03 and the Grant
for Support of Scientific
Schools NSh-8065.2006.2.

\subsection{Notation}
Our underlying symmetric monoidal categories
are the category of graded vector spaces
and the category of chain complexes.
Objects of these categories we
sometimes loosely call ``spaces''.
By suspension $\Si V$ of a space $V$
we mean $\ve \otimes V$, where $\ve$ is a
one-dimensional vector space placed in degree
$+1$\,. The underlying field $\bbK$
has characteristic zero.

For an operad $\cO$
we denote by $\Alg_{\cO}$ the category of algebras
over the operad $\cO$. Dually, for a cooperad
$\cC$ we denote by $\Coalg_{\cC}$ the category
of nilpotent\footnote{For the definition of
{\it nilpotent} coalgebra see section $2.4.1$ in
\cite{Dima11}.}
coalgebras over the cooperad $\cC$\,.
By {\it corestriction} we mean the canonical
map
\begin{equation}
\label{corest}
u_V : \bbF_{\cC}(V) \to V
\end{equation}
from the cofree coalgebra
$\bbF_{\cC}(V)$ to the space of its
cogenerators $V$\,.
$Ho(\cO)$
is reserved for the operad of homotopy $\cO$-algebras.

For a quadratic operad $\cO$ we will denote by
$\cO^{\vee}$ the Koszul dual cooperad \cite{Fresse},
\cite{GK}.
$Bar(\cO)$ will denote the
bar construction of the (augmented) operad $\cO$,
and $Cobar(\cC)$ will denote the cobar construction
of the (coaugmented) cooperad $\cC$\,.
For an operad $\cO$ (resp. cooperad $\cC$) and
a space $V$ we denote by
$\bbF_{\cO}(V)$ (resp. by $\bbF_{\cC}(V)$) the
free algebra (resp. cofree coalgebra) over the
operad $\cO$ (resp. cooperad $\cC$).
We will use the following list of
(co)operads:
\begin{itemize}

\item $\ass$ (resp. $\coass$) is the
operad of associative algebras
\underline{without unit} (resp. the
cooperad of coassociative algebras
without counit),

\item $\lie$ (resp. $\colie$) is
the operad of Lie algebras (resp. the cooperad
of Lie coalgebras),

\item $\comm$ (resp. $\cocomm$)
is the operad of commutative (associative)
algebras (resp. the
operad of cocommutative coassociative
coalgebras),

\item $\Ger$ denotes the operad of Gerstenhaber
algebras,

\item $\Bra$ denotes the operad of
braces \cite{Ezra}, \cite{Gruzia}.

\item $\lie^+$ (resp. $\colie^+$) denotes the
$2$-colored operad of pairs ``Lie algebra $+$ its
module'' (resp. the $2$-colored cooperad of pairs
``Lie coalgebra $+$ its comodule'')

\item $\comm^+$ (resp. $\cocomm^+$) denotes the
$2$-colored operad of pairs ``commutative algebra $+$ its
module'' (resp. the $2$-colored cooperad of pairs
``cocommutative coalgebra $+$ its comodule'')

\end{itemize}

By ``suspension'' of a (co)operad $\cO$
we mean the (co)operad $\La(\cO)$
whose $m$-th space is
\begin{equation}
\label{susp-op}
\La(\cO)(m) = \Si^{1-m} \cO(m) \otimes \sgn_{m}\,,
\end{equation}
where $\sgn_{m}$ is the sign representation of
the symmetric group $S_m$\,.

For a commutative algebra $\cB$ and a $\cB$-module
$\cV$ we denote by $S_{\cB}(\cV)$ the symmetric algebra
of $\cV$ over $\cB$\,.
$S_{\cB}^m(\cV)$ stands for the $m$-th component
of this algebra. If $\cB=\bbK$ then $\cB$ is omitted
from the notation.

For an associative algebra $A$
$$
\Cbu(A) = Hom (A^{\otimes \bul}, A)
$$
denotes the normalized Hochschild cochain complex.
$[ , ]_G$ and $\cdot$ denote the Gerstenhaber
bracket and the cup product in $\Cbu(A)$,
respectively.

\subsection{Preliminaries}
In this section we recall some basic
facts about (co)algebras and (co)operads
we will use in the paper.

\subsubsection{Koszul (co)operads and their
(co)algebras}
Let us recall \cite{GK} that
a quadratic operad
$\cO$ is {\it Koszul}
if the canonical map
\begin{equation}
\label{izcobar}
Cobar(\cO^{\vee}) \to \cO
\end{equation}
is a quasi-isomorphism.

For any coaugmented
cooperad $\cC$ the operad $Cobar(\cC)$ is obviously
cofibrant \cite{Hinich}. Thus, as well as Koszul
algebras \cite{Priddy}, Koszul operads
\cite{GK} admit simpler cofibrant
resolutions.

This observation motivates
the following definition:
\begin{defi}
\label{O-inf}
For a Koszul operad $\cO$
the operad $Cobar(\cO^{\vee})$ can
be viewed as an operad $Ho(\cO)$ of homotopy
$\cO$-algebras.
\end{defi}

Let us also recall that for any
coaugmented cooperad $\cC$ the algebras
over the operad $Cobar(\cC)$ have the
following explicit description
due to \cite{GJ}:
\begin{pred}[Proposition 2.15, \cite{GJ}]
\label{GJ}
Let $\cC$ be a coaugmented cooperad
and let
$V$ be a space. There is a natural
bijection between the degree $1$ codifferentials
on the cofree coalgebra $\bbF_{\cC}(V)$ and
the $Cobar(\cC)$-algebra structures
on $V$.
\end{pred}

Furthermore, one can describe coderivations
of a cofree coalgebra
\begin{pred}[Proposition 2.14 \cite{GJ}]
\label{coder-cofree}
Let $\cC$ be a cooperad and let
$Coder(\bbF_{\cC}(V))$ denote the Lie
algebra of coderivations of
the cofree coalgebra $\bbF_{\cC}(V)$.
Then the composition with the corestriction
$u_V: \bbF_{\cC}(V) \to V$ induces an
isomorphism
\begin{equation}
\label{der-free}
Coder(\bbF_{\cC}(V)) \cong
\Hom(\bbF_{\cC}(V), V)
\end{equation}
of the graded vector space of coderivations
of the cofree coalgebra $\bbF_{\cC}(V)$
with the graded vector space
$\Hom(\bbF_{\cC}(V), V)$\,.
\end{pred}

In view of definition \ref{O-inf}
proposition \ref{GJ} implies that
\begin{cor}
\label{O-inf1}
For a Koszul operad $\cO$ and
a space $V$ there is a natural bijection
between the $Ho(\cO)$-algebra structures
on $V$ and
degree $1$ codifferentials on the
cofree coalgebra $\bbF_{\cO^{\vee}}(V)$\,.
\end{cor}

The above corollary motivates the
definition of homotopy morphisms between 
$Ho(\cO)$-algebras for a Koszul operad $\cO$:
\begin{defi}
\label{def-morph}
A homotopy morphism $F$ from a
$Ho(\cO)$-algebra $V_1$ to a $Ho(\cO)$-algebra
$V_2$ is a morphism of the
DG coalgebras
\begin{equation}
\label{morph-F}
F : (\bbF_{\cO^{\vee}}(V_1), d_{V_1}) \to
(\bbF_{\cO^{\vee}}(V_2), d_{V_2})\,,
\end{equation}
where $d_{V_1}$ and $d_{V_2}$ are codifferentials
corresponding to the $Ho(\cO)$-algebra structures
on $V_1$ and $V_2$, respectively.
\end{defi}
We reserve the following notation for such morphisms
\begin{equation}
\label{notation}
F : V_1 \brarrow V_2\,.
\end{equation}
Let $u$ and $\la$ be, respectively,
the corestriction (\ref{corest}) and
the coaugmentation of
the cooperad $\cO^{\vee}$\,.
Similarly to proposition \ref{coder-cofree}
one can show that every morphism of
cofree coalgebras
$$
F : \bbF_{\cO^{\vee}}(V_1)  \to
\bbF_{\cO^{\vee}}(V_2)
$$
is uniquely determined by its composition
with the corestriction
$$
u_{V_2} : \bbF_{\cO^{\vee}}(V_2) \to V_2\,.
$$
In what follows we call the components
\begin{equation}
\label{str-maps}
F_n = u_{V_2} \circ F \Big|_{\cO^{\vee}(n) \otimes_{S_n} V_1^{\otimes\, n}}
\cO^{\vee}(n) \otimes_{S_n} V_1^{\otimes\, n} \to V_2
\end{equation}
of the composition $u_{V_2} \circ F$ {\it the structure maps}
of the homotopy morphism $F$ (\ref{morph-F}).

Using definition \ref{def-morph} it is
not hard to show that the first structure map
(\ref{str-maps}) of (\ref{morph-F})
\begin{equation}
\label{F-1}
F_1 = u_{V_2} \circ F \circ \la : V_1 \to V_2
\end{equation}
is always a morphism of the corresponding complexes.
This observation motivates the following
definition:
\begin{defi}
\label{q-iso}
Let $\cO$ be a Koszul operad.
A quasi-isomorphism between $Ho(\cO)$-algebras
$V_1$ and $V_2$ is a homotopy morphism
(\ref{morph-F}) whose first structure map
(\ref{F-1}) is a quasi-isomorphism of
complexes.
\end{defi}

We will often have to deal with
(co)free (co)algebras over a suspended
(co)operad. The following easy
fact will be very useful:
\begin{equation}
\label{La-Si}
\bbF_{\La\cO}(V) = \Si\, \bbF_{\cO} (\Si^{-1} V)\,.
\end{equation}
Equation (\ref{La-Si}) holds both for
operads and cooperads.

Let $\cO$ be an augmented operad and let
$\cO^{\vee}$ be its Koszul dual cooperad.
Then according to \cite{GJ} we have
\begin{teo}[Theorem 2.25, \cite{GJ}]
\label{adjun1}
There is an adjunction
\begin{equation}
\label{adjun}
\Om_{\cO} \, :\, \Coalg_{\cO^{\vee}} \rightleftharpoons
\Alg_{\cO} \, : \, B_{\cO^{\vee}}
\end{equation}
between the category $\Coalg_{\cO^{\vee}}$
of nilpotent $\cO^{\vee}$-coalgebras and
the category $\Alg_{\cO}$ of $\cO$-algebras.
If $\cO$ is Koszul the counit
\begin{equation}
\label{counit1}
\eta_{\cO} : \Om_{\cO} (B_{\cO^{\vee}}(A)) \to A
\end{equation}
and the unit
\begin{equation}
\label{unit1}
\ve_{\cO^{\vee}} : C \to  B_{\cO^{\vee}} ( \Om_{\cO}(C))
\end{equation}
of this adjunction are quasi-isomorphisms.
\end{teo}
{\bf Remark 1.} Notice that, if $A$ is an $\cO$-algebra
then as a coalgebra $B_{\cO^{\vee}}(A)$ is just
the cofree coalgebra $\bbF_{\cO^{\vee}}(A)$
cogenerated by $A$\,. Dually, if $C$ is
a $\cO^{\vee}$-coalgebra then as an algebra
$\Om_{\cO}(C)$
is $\bbF_{\cO}(C)$\,. Thus the algebra
structure on $A$ (resp. the coalgebra structure
on $C$) contributes only to the codifferential
on $B_{\cO^{\vee}}(A)$ (resp. differential
on $\Om_{\cO}(C)$)

~\\
{\bf Remark 2.} An obvious analogue of theorem
\ref{adjun1} can be proved for colored
operads.

The counit (\ref{counit1}) and
the unit (\ref{unit1}) of the adjunction
(\ref{adjun}) admit simple explicit
descriptions\footnote{See sections $2.3$ and
$2.4$ in \cite{GJ}.}. To describe the map
$\eta_{\cO}$ (\ref{counit1}) we notice that
it is a map of $\cO$-algebras
from a free $\cO$-algebra. Hence, $\eta_{\cO}$
is uniquely determined by its restriction to the
space of generators. Using the definition
of the functors $\Om_{\cO}$ and $B_{\cO^{\vee}}$
 \cite{GJ}
it is not hard to see that
\begin{equation}
\label{def-eta}
\eta_{\cO} \Big|_{B_{\cO^{\vee}}(A)}  = u_A\,,
\end{equation}
where $u_A$ is the corestriction
$$
u_A : \bbF_{\cO^{\vee}} (A) \to A\,.
$$
Similarly, the map $\ve_{\cO^{\vee}}$
being a map to a cofree $\cO^{\vee}$-coalgebra
is uniquely determined by its composition
$u_{\bbF_{\cO}(C)} \circ \ve_{\cO^{\vee}}$
with the corestriction
$$
u_{\bbF_{\cO}(C)} : \bbF_{\cO^{\vee}}(\bbF_{\cO}(C)) \to
\bbF_{\cO}(C)
$$
onto cogenerators.

Again, using the definition
of the functors $\Om_{\cO}$ and $B_{\cO^{\vee}}$
 \cite{GJ}
one can show that
\begin{equation}
\label{def-ve}
u_{\bbF_{\cO}(C)} \circ \ve_{\cO^{\vee}} = \rho\,,
\end{equation}
where $\rho$ is the embedding of
generators to the free algebra:
$$
\rho: C \to \bbF_{\cO}(C)\,.
$$

Let us now prove that if
$\cO$ is a Koszul operad then for
every $Ho(\cO)$-algebra $V$ one
can construct a DG $\cO$-algebra which
is quasi-isomorphic\footnote{An analogous
construction in topology is called
{\it rectification} \cite{BVogt}.}
to the initial
$Ho(\cO)$-algebra $V$.

To construct this algebra we first notice
that, due to corollary \ref{O-inf1},
a $Ho(\cO)$-algebra $V$ corresponds to
a degree $1$ codifferential $D$ of the
cofree  $\cO^{\vee}$-coalgebra
$\bbF_{\cO^{\vee}}(V)$. Thus a $Ho(\cO)$-algebra $V$
gives us the DG $\cO^{\vee}$-coalgebra
$$
(\bbF_{\cO^{\vee}}(V), D)\,.
$$

Applying theorem \ref{adjun1} to the
above coalgebra we get the
following morphism of DG coalgebras
\begin{equation}
\label{ve-CV}
\chi^V =
\ve_{\cO^{\vee}} : (\bbF_{\cO^{\vee}}(V), D)
\to B_{\cO^{\vee}} (\Om_{\cO}(  \bbF_{\cO^{\vee}}(V), D ))\,.
\end{equation}
Definition \ref{def-morph}
allows us to treat (\ref{ve-CV}) as a
homotopy morphism from $V$ to
$\Om_{\cO}( \bbF_{\cO^{\vee}}(V), D )$\,.
The following proposition shows that this homotopy
morphism is a quasi-isomorphism
in the sense of definition \ref{q-iso}.
\begin{pred}
\label{pro-chi}
Let $\cO$ be a Koszul operad and
$D$ be a codifferential on the
cofree $\cO^{\vee}$-coalgebra
$\bbF_{\cO^{\vee}}(V)$ defining a
$Ho(\cO)$-algebra structure on a space $V$.
Then the homotopy morphism $\chi^V$
\begin{equation}
\label{chi-V}
\chi^V : V \brarrow
 \Om_{\cO}(\bbF_{\cO^{\vee}}(V), D)\,.
\end{equation}
corresponding to the map (\ref{ve-CV})
is a quasi-isomorphism.
\end{pred}
{\bf Proof.} Let, as above,
$u$ and $\la$ denote, respectively,
the corestriction (\ref{corest}) and the coaugmentation of
the cooperad $\cO^{\vee}$\,.

According to definition \ref{q-iso}
we have to prove that the composition
$$
\chi^V_1 = u \circ \chi^V \circ \la
: V \to \Om_{\cO}( \bbF_{\cO^{\vee}}(V), D )
$$
is a quasi-isomorphism of the corresponding
complexes, where the differential
on $V$ is the composition:
\begin{equation}
\label{D-1}
D_1 = u \circ D \circ \la\,.
\end{equation}

Due to equation (\ref{def-ve}) the composition
$u \circ \chi_V$ is simply the embedding
$$
\bbF_{\cO^{\vee}}(V) \hookrightarrow
\bbF_{\cO} (\bbF_{\cO^{\vee}}(V))
$$
of the space $\bbF_{\cO^{\vee}}(V)$ of
generators to the free $\cO$-algebra
$\bbF_{\cO} (\bbF_{\cO^{\vee}}(V))$\,.
Hence, $\chi^V_1$ coincides with the embedding
of $V$ into $\Om_{\cO}(\bbF_{\cO^{\vee}}(V), D)$\,.
In particular, the graded vector space $V$
with the differential $D_1$ (\ref{D-1}) is
a subcomplex of $\Om_{\cO}(\bbF_{\cO^{\vee}}(V), D)$\,.

It remains to prove that the quotient complex
of the complex $\Om_{\cO}(\bbF_{\cO^{\vee}}(V), D)$
by the subcomplex $(V, D_1)$ is acyclic.
Using the filtration corresponding to the total
degree of the space $\bbF_{\cO^{\vee}}(V)$
it is not hard to reduce
the latter question to the computation of the
homology groups of the cofree $\cO^{\vee}$-coalgebra
$\bbF_{\cO^{\vee}}(V)$. Then the desired statement
follows easily from the criterion of Ginzburg and
Kapranov \cite{GK} (theorem $4.2.5$) and the koszulity of
the operad $\cO$\,. $\Box$

\subsubsection{Operads $\lie$, $\comm$, $\Ger$,
et cetera} \label{razdel}
Let us recall the following
well known facts
\begin{pred}[V. Ginzburg and M. Kapranov,
\cite{GK}]
\label{lie-comm}
The operads $\lie$ and $\comm$ are Koszul.
\begin{equation}
\lie^{\vee} = \La (\cocomm)\,, \qquad
\comm^{\vee} = \La (\colie)\,.
\end{equation}
\end{pred}
Furthermore,
\begin{pred}[E. Getzler and J.D.S Jones,
\cite{GJ}]\label{Ger-Koszul}
The operad $\Ger$ is Koszul.
\end{pred}
{\bf Remark.} In \cite{GJ} E. Getzler and
J.D.S. Jones proved koszulity for wider
class of operads which includes the
operad $\Ger$ of Gerstenhaber algebras.
In fact, using the Ginzburg-Kapranov
criterion of koszulity (see \cite{GK},
Theorem $4.2.5$) it is not hard to
prove proposition \ref{Ger-Koszul}.
This proof can be found in
\cite{Dima11} (see subsection $5.4.6$).

A simple computation shows that
coalgebras over $\Ger^{\vee}$ are
equipped with a cocommutative coassociative
coproduct of degree $-2$ and a
Lie cobracket of degree $-1$, satisfying
the Leibniz rule. Thus
\begin{equation}
\label{free-coGer}
\bbF_{\Ger^{\vee}}(V) = \Si^{2}
\bbF_{\cocomm} (\Si^{-1} \bbF_{\colie}(\Si^{-1}V))\,,
\end{equation}
and thanks to corollary \ref{O-inf1} and
proposition \ref{Ger-Koszul} the homotopy
Gerstenhaber algebra structures on
a space $V$ are in a one-to-one
correspondence with the degree $1$ codifferentials
of the coalgebra (\ref{free-coGer})\,.

Let $(V,\, [\,,\,]_V\,,\, \bul_V)$ be a 
Gerstenhaber algebra with the bracket 
$[\,,\,]_V$ and the product $\bul_V$\,. 
Using the fact that $V$ is both a commutative 
algebra and a $\La\lie$-algebra one can 
split the construction of the
DG $\Ger^{\vee}$-coalgebra
$B_{\Ger^{\vee}}(V)$ into two steps. 
 
In the first step we upgrade the 
DG $\La\colie$-coalgebra $B_{\La\colie}(V,\, \bul_V)$
to a DG Lie bialgebra (with the existing
cobracket of degree $1$ a Lie bracket 
$[\,,\,]$ of degree $-1$) using the 
bracket $[\,,\,]_V$ on $V$\,.
The bracket $[\,,\,]$ on $B_{\La\colie}(V,\, \bul_V)$ 
is uniquely determined by the formula
\begin{equation}
\label{Dima}
u_V([ X , Y ]) = [u_V(X) , u_V(Y)]_V\,, 
\qquad X, Y \in B_{\La\colie}(V,\, \bul_V)
\end{equation}
and the compatibility condition 
with the existing cobracket. Here, as 
above, $u_V$ denotes the corestriction 
$u_V : B_{\La\colie}(V) \to V$ on the 
space $V$ of cogenerators. The resulting 
bracket $[\,,\,]$ is compatible with the 
codifferential on  $B_{\La\colie}(V,\, \bul_V)$ 
due to the Leibniz rule for the bracket $[\,,\,]_V$
on $V$.

In the second step we apply the functor 
$B_{\La^2\cocomm}$ to the Lie bialgebra 
$B_{\La\colie}(V,\, \bul_V)$, regarding 
the latter as a DG $\La\lie$-algebra, and 
obtain a DG $\La^2\cocomm$-coalgebra 
\begin{equation}
\label{Dima1}
C = B_{\La^2\cocomm}\, B_{\La\colie}(V,\, \bul_V)\,.
\end{equation}
Then extending in the similar way 
the Lie cobracket from $B_{\La\colie}(V,\, \bul_V)$
to $C$ using the compatibility with
the coproduct on $C$ we get the desired 
DG $\Ger^{\vee}$-coalgebra 
\begin{equation}
\label{Dima11}
B_{\Ger^{\vee}}(V) = 
B_{\La^2\cocomm}\, B_{\La\colie}(V,\, \bul_V)\,.
\end{equation}
The compatibility of 
the resulting cobracket on $C$ (\ref{Dima1}) with 
the codifferential
follows from the compatibility of the 
cobracket and the bracket
on $B_{\La\colie}(V,\, \bul_V)$\,.

The above contruction of the 
DG $\Ger^{\vee}$-coalgebra (\ref{Dima11})
shows that the functors $B_{\La\colie}$
and $B_{\La^2\cocomm}$
can be extended to functors 
\begin{equation}
\label{B'}
B'_{\La\colie} : \Alg_{\Ger}
\to \Bialg
\end{equation}
and 
\begin{equation}
\label{B''}
B'_{\La^2\cocomm} : \Bialg \to 
\Coalg_{\Ger^{\vee}}\,,
\end{equation}
where $\Bialg$ is the category of DG Lie bialgebras 
with the bracket of degree $-1$ and the nilpotent
cobracket of degree $1$. The composition of 
these two functors gives us the 
functor $B_{\Ger^{\vee}}$
\begin{equation}
\label{Dima111}
B_{\Ger^{\vee}} =
B'_{\La^2\cocomm}\circ B'_{\La\colie}\,.
\end{equation}
In what follows we will omit the prime 
in the notations of the functors (\ref{B'})
and (\ref{B''})\,.

It is also clear from (\ref{free-coGer}) that
there are canonical maps of operads:
\begin{equation}
\label{lie-comm-Ger}
Ho(\La(\lie)) \to \Hoger\,, \qquad
\Hocomm \to \Hoger\,.
\end{equation}

We use the following realization of the cofree
Lie coalgebra
\begin{equation}
\label{Colie}
\bbF_{\colie}(V) =
\bbF_{\coass} (V)\, / \,
\bbF_{\coass} (V) \, \bul_{s h}\, \bbF_{\coass} (V)\,,
\end{equation}
where $\bul_{s h}$ denotes the shuffle product
\begin{equation}
\label{shuffle}
\lan a_1, a_2, \dots, a_p \ran
\, \bul_{s h} \,
\lan a_{p+1}, a_{p+2}, \dots, a_{p+q} \ran =
\sum_{\ve \in  S h (p,q)}
\lan a_{\ve^{-1}(1)}, a_{\ve^{-1}(2)},
\dots
a_{\ve^{-1}(p+q)} \ran \,,
\end{equation}
$$
a_i \in  V\,.
$$

Let us conclude this section with the following statement
\begin{pred}
\label{comm-Koszul}
The operad $\comm^+$
of pairs ``commutative algebra $+$ its module''
is Koszul.
\end{pred}
{\bf Proof.} It is not hard to see that
\begin{equation}
\label{comm-Kos}
(\comm^+)^{\vee} = \La \colie^+\,.
\end{equation}
Hence, due to a proper version of
Ginzburg-Kapranov criterion (theorem $4.2.5$ in \cite{GK})
it suffices to show that Harrison homology groups
of the free commutative algebra with trivial
coefficients and with coefficients in itself
are concentrated in the lowest degrees.
Due to the realization (\ref{Colie})
the latter question reduces to
a simple computation with Hochschild chain complexes
of the free commutative algebra. $\Box$

\section{The $\Hoger$-algebra structure
on the Hochschild complex}
It follows from the result of
\cite{M-Smith} and the formality of the
operad of little squares \cite{Dima-Disk}
that we have a quasi-isomorphism
\begin{equation}
\label{vesch}
\Hoger \erarrow \Bra
\end{equation}
from the operad $\Hoger$ of
homotopy Gerstenhaber
algebras to the operad of braces
$\Bra$\,.

This result implies that the normalized Hochschild
complex $C^{\bul}(A)$ of any associative
algebra $A$ is naturally a $\Hoger$-algebra.
In what follows we say {\it the}
$\Hoger$-algebra structure on $C^{\bul}(A)$
referring to the structure provided by
the map (\ref{vesch}).
In this section we discuss properties of
this $\Hoger$-algebra.

The first remarkable property of the
homotopy Gerstenhaber algebra $C^{\bul}(A)$
can be formulated as
\begin{teo}[\cite{Dima1}]
\label{non-spoiled}
The $Ho(\La(\lie))$-algebra structure on
$C^{\bul}(A)$ induced by the first map
in (\ref{lie-comm-Ger}) coincides with
the DG $\La(\lie)$-algebra structure
given by the Hochschild differential
and the Gerstenhaber bracket.
\end{teo}
{\bf Proof.} If $\cO$ is a Koszul operad
then so is $\La \cO$ and $\La^{-1}\cO$.
Furthermore, it is easy to see that
the suspension commutes with the operation
of taking the Koszul dual (co)operad:
$$
(\La(\cO))^{\vee}= (\La(\cO)^{\vee})\,.
$$

Thus the $Ho(\La(\lie))$-algebra structure
on $\Cbu(A)$ is determined by a
degree $1$ map
$$
Q: \bbF_{\La^2\cocomm}(\Cbu(A)) \to \Cbu(A)
$$
or, in other words, by an infinite collection
of maps
$$
Q_m : S^m (C^{\bul}(A)) \to C^{\bul}(A)
$$
of degree $3-2m$\,.

To prove the theorem we need to show
that
\begin{equation}
\label{Q-m}
Q_m=0\,, \qquad  m > 2
\end{equation}
and
\begin{equation}
\label{Q-2=G}
Q_2 (P_1, P_2) = (-)^{|P_1|} [P_1, P_2]_G\,,
\qquad \forall~ P_i \in \Cbu(A)\,,
\end{equation}
where $[,]_G$ is the Gerstenhaber bracket
\cite{Ger}\,.

To prove (\ref{Q-m})
we notice that the lowest degree
$m$-ary operation
we can get by combining the braces
and the cup product is
\begin{equation}
\label{lowest}
(P_1, \dots, P_m) \to P_i \{P_1, \dots, \hat{P_i}
\dots, P_m \}\,,
\end{equation}
where $i$ runs from $1$ to $m$ and
$P_i\in C^{\bul}(A)$\,.

It is easy to see that the degree of
the operation (\ref{lowest}) is $1-m$.
Thus if $m>2$ then the degree of $Q_m$
is lower than the degree that can be obtained
by combining the braces and the cup-product.
Hence for any $m>2$, $Q_m=0$.

Since on the cohomology space $HH^{\bul}(A)$
both $\Hoger$-algebras induce the same
graded Lie algebra structure we conclude
that for any pair $P_1, P_2 \in \Cbu(A)$
\begin{equation}
\label{Q-2=G1}
Q_2 (P_1, P_2) = (-)^{|P_1|} [P_1, P_2]_G +
\pa \Psi(P_1, P_2) - \Psi(\pa P_1, P_2)
- (-)^{|P_1|} \Psi(P_1, \pa P_2)\,,
\end{equation}
where $\pa$ is the Hochschild coboundary
operator and
$$
\Psi : \Cbu(A)\otimes \Cbu(A) \to \Cbu(A)
$$
is a map of degree $-2$ expressed in terms of
braces and the cup-product.

Using the braces and the cup-product it
is impossible to get a binary operation of
degree $-2$\,. This observation implies
equation (\ref{Q-2=G}) and concludes
the proof. $\Box$

Due to corollary \ref{O-inf1} and
proposition \ref{Ger-Koszul} the $\Hoger$-algebra
on $C^{\bul}(A)$ is encoded by a
degree $1$ codifferential
\begin{equation}
\label{M}
M : \bbF_{\Ger^{\vee}}(C^{\bul}(A))
\to  \bbF_{\Ger^{\vee}}(C^{\bul}(A))
\end{equation}
of the cofree  $\Ger^{\vee}$-coalgebra
$\bbF_{\Ger^{\vee}}(C^{\bul}(A))$
cogenerated by $\Cbu(A)$\,.

Thanks to proposition \ref{coder-cofree}
$M$ is uniquely determined by the
composition $ u_{\Cbu(A)} \circ M$ with the
corestriction
\begin{equation}
\label{counit}
u_{\Cbu(A)}  : \bbF_{\Ger^{\vee}}(C^{\bul}(A)) \to \Cbu(A)\,.
\end{equation}
We denote this composition $ u_{\Cbu(A)} \circ M $
by $m$
\begin{equation}
\label{m}
m = u_{\Cbu(A)} \circ M : \Si^2 \bbF_{\cocomm}(\Si^{-1} \bbF_{\colie}(\Si^{-1}\Cbu(A)))
\to \Cbu(A)\,.
\end{equation}
As well as (\ref{M}) the map (\ref{m})
is of degree $1$\,.

Let $\Der(A)$ denote the Lie algebra of derivations
of $A$ and let
\begin{equation}
\label{polyvec}
V^{\bul}(A)= S^{\bul}_A(\Si\,\Der(A) )\,.
\end{equation}
be the Gerstenhaber algebra of polyvectors
on $Spec(A)$\,.

It is well-known \cite{HKR} that the natural
embedding
\begin{equation}
\label{HKR}
V^{\bul}(A) \to \Cbu(A)
\end{equation}
is a quasi-isomorphism  of complexes
where $V^{\bul}(A)$ is viewed as a complex
with the vanishing differential.
In other words, $V^{\bul}(A)$ is the cohomology
space of the Hochschild cochain complex
$(\Cbu(A), \pa)$. In \cite{HKR} (see theorem $5.2$)
it was shown that the product induced on
the cohomology space $V^{\bul}(A)$ coincides with
natural product on (\ref{polyvec}).
Furthermore, it is not hard to show that
the $\La\lie$-algebra structure induced on the cohomology
space $V^{\bul}(A)$ coincides with the one given
by the Schouten-Nijenhuis bracket \cite{Nijenhuis}.

Let us introduce the following
sub $\Ger^{\vee}$-coalgebra
$$
\Xi(A) \subset \bbF_{\Ger^{\vee}}(\Cbu(A))
$$
\begin{equation}
\label{DT}
\Xi(A) = \bbF_{\La^{2}\cocomm} \circ
\bbF_{\La\colie^+} (A, \Si\Der(A))\,,
\end{equation}
where
$$
\bbF_{\La\colie^+} (A, \Si\Der(A)) =
\bbF_{\La\colie}(A)
\oplus \bbF_{\La\colie^+} (A, \Si\Der(A))^+
$$
is the direct sum of the cofree
$\La\colie$-coalgebra $\bbF_{\La\colie}(A)$
cogenerated by $A$
and the cofree comodule
$\bbF_{\La\colie^+} (A, \Si\Der(A))^+$
over $\bbF_{\La\colie}(A)$ cogenerated by
$\Der(A)$ placed in degree $1$\,.

It is not hard to see that
$\Xi(A)$ consists of sums of
expressions:
\begin{equation}
\label{DT1}
(\la_1 \otimes p_1)^{S_{|p_1|}}
\dots
(\la_k \otimes p_k)^{S_{|p_k|}}
\end{equation}
where $p_i$ are monomials in
$\bigotimes(A\oplus \Der(A))$ in which
an element of $\Der(A)$ can appear at most once,
$\la_i$ is an element of $\La\colie(|p_i|)$
and $|p_i|$ denotes the degree of $p_i$\,.

We conclude this section with the
following theorem:
\begin{teo}
\label{ona}
The map $m$ (\ref{m}) vanishes on monomials
of $\Xi(A)$ which have more than two
components. Furthermore,
\begin{equation}
\label{binary}
m (P_1 \, P_2) = (-)^{|P_1|}[P_1, P_2]_G\,,
\qquad
m \lan P_1 , P_2 \ran = (-)^{|P_1|} P_1 \cdot P_2
\end{equation}
whenever the monomial
$P_1 \, P_2$ (resp. the monomial
$\lan P_1 , P_2 \ran$ ) belongs to $\Xi(A)$\,.
\end{teo}
{\bf Proof.}
Equation (\ref{binary})  follows easily
from the degree bookkeeping and it
remains to prove that the
map $m$ (\ref{m}) vanishes on monomials
of $\Xi(A)$ which have more than two
components. The latter is equivalent to the
following equations
\begin{equation}
\label{vff}
m  \lan v, f_1, f_2 \ran =0\,,
\qquad
m \lan f_1, v, f_2 \ran =
m \lan f_1, f_2, v \ran =0\,,
\end{equation}
\begin{equation}
\label{vfv}
m  ( \lan v_1, f \ran \, v_2 ) =0\,,
\end{equation}
\begin{equation}
\label{vvv}
m  (v v_1 v_2)=0\,,
\end{equation}
$$
\forall \qquad
v, v_1, v_2 \in \Der(A) \qquad
f, f_1, f_2 \in A\,,
$$
where we use the notation
$$
\lan P_1, \dots, P_m \ran\,, \qquad
P_i \in \Cbu(A)
$$
for elements in
$\bbF_{\colie}(\Si^{-1} \Cbu(A))$
keeping in mind the realization
(\ref{Colie}).

Indeed, any other higher combination will
have degree $<-1$ while the complex $\Cbu(A)$ lies
in non-negative degrees\,. Since the map
(\ref{m}) is of degree one all these
combinations get sent to zero.

Equation (\ref{vvv}) follows from theorem
\ref{non-spoiled} and
the second pair of equations in (\ref{vff})
follows from the first equation.
The latter can be proved using the
fact that $m$ has to vanish on
shuffle-products (\ref{shuffle}).

Thus we are left with the pair
of the following equations:
\begin{equation}
\label{vff-vfv}
m  \lan v, f_1, f_2 \ran =0\,,
\qquad
m  \lan v_1, f \ran \, v_2 =0\,.
\end{equation}

Since the map $m$ (\ref{m}) is of degree $1$
and $\lan v, f_1, f_2 \ran$ is of degree
$-1$ in the coalgebra
$$
\Si^2 \bbF_{\cocomm}(\Si^{-1}\bbF_{\colie}
(\Si^{-1} \Cbu(A)))
$$
the element
$$
m  \lan v, f_1, f_2 \ran \in C^0(A)=A\,.
$$

Thus, using the fact that the map $m$ (\ref{m}) is
given in terms of the brace operations
\cite{Ezra}, \cite{Gruzia}, we conclude that
the most general expression for
$m  \lan v, f_1, f_2 \ran$ is
\begin{equation}
\label{m-vff}
m  \lan v, f_1, f_2 \ran =
\al f_1 v(f_2) + \beta f_2 v(f_1)\,,
\end{equation}
where $\al, \beta\in \bbK$ and
$v(\cdot)$ denotes the action of the
derivation on a function.

On the other hand, applying the corestriction
$ u_{\Cbu(A)} $ (\ref{counit}) to the
equation
$$
M^2 (\lan v, f_1, f_2, f_3  \ran ) =0
$$
and using (\ref{binary}) we get
\begin{equation}
\label{dushki}
f_3 \cdot m \lan v, f_1, f_2 \ran -
 m \lan f_1 \cdot v, f_2, f_3 \ran
+  m \lan v, f_1 \cdot f_2, f_3 \ran
-  m \lan v, f_1, f_2 \cdot f_3 \ran =0\,,
\end{equation}
where $\cdot$ denotes the cup product in
$\Cbu(A)$\,.

Due to (\ref{m-vff}) equation (\ref{dushki})
boils down to
\begin{equation}
\label{al-beta}
\beta \, f_2\cdot f_3 \cdot v(f_1) -
\al \, f_1 \cdot f_2 \cdot v(f_3) =0\,.
\end{equation}
Since (\ref{al-beta}) holds for any
associative algebra $A$ the coefficients
$\al$ and $\beta$ vanish and the first
equation in (\ref{vff-vfv}) is satisfied.

Similar degree bookkeeping shows that
$$
m ( \lan v_1, f \ran \,  v_2 ) \in C^0(A)=A\,.
$$
Hence, the most general expression for
$ m ( \lan v_1, f\ran\, v_2 ) $ is
\begin{equation}
\label{m-vfv}
m  ( \lan v_1, f\ran\, v_2 )  =
\mu v_2 (v_1 (f)) + \nu v_1 (v_2(f))\,,
\end{equation}
where $\mu, \nu\in \bbK$\,.

Applying the corestriction $u_{\Cbu(A)}$ (\ref{counit}) to
the equation
$$
M^2 (\lan v_1, f \ran\, v_2 \, v_3 ) =0
$$
and using (\ref{binary}) we get
\begin{equation}
\label{dushki1}
v_2 (m  ( \lan v_1, f\ran\, v_3 ) ) +
m  ( \lan [v_1, v_2]_G, f\ran\, v_3 )
+ m  ( \lan v_1, v_3(f) \ran\, v_2 ) -
(2 \leftrightarrow 3) = 0\,.
\end{equation}

Due to (\ref{m-vfv}) equation (\ref{dushki1})
boils down to
\begin{equation}
\label{mu-nu}
\mu \, [v_2, v_3]_G (v_1(f)) +
\nu \, v_1 ([v_2, v_3]_G (f)) =0\,.
\end{equation}
Since equation (\ref{mu-nu}) holds
for any associative algebra $A$
the coefficients $\mu$ and $\nu$
necessarily vanish and the second equation
in (\ref{vff-vfv}) is satisfied.

The theorem is proved.  $\Box$

\section{The formality theorem}
In this section we construct a chain of
quasi-isomorphisms between the homotopy Gerstenhaber
algebra $\Cbu(A)$ of Hochschild cochains
of a regular commutative algebra $A$ and
the Gerstenhaber algebra $V^{\bul}(A)$ (\ref{polyvec})
of multiderivations of $A$\,.

First, we observe that, as in
equation (\ref{ve-CV}), the map (\ref{unit1})
gives us a homotopy morphism
\begin{equation}
\label{chi-C}
\chi^{\Cbu(A)} :
\Cbu(A) \brarrow \Om_{\Ger} (\bbF_{\Ger^{\vee}}(\Cbu(A)), M)\,,
\end{equation}
where $M$ is the codifferential of the
coalgebra $\bbF_{\Ger^{\vee}}(\Cbu(A))$ defining the
homotopy Gerstenhaber algebra structure on $\Cbu(A)$\,.
Due to proposition \ref{pro-chi} the homotopy morphism
(\ref{chi-C}) is a quasi-isomorphism.

Second, theorem \ref{adjun1} gives a quasi-isomorphism
of DG Gerstenhaber algebras
\begin{equation}
\label{eta-V}
\eta_{\Ger} :
\Om_{\Ger} (B_{\Ger^{\vee}}( V^{\bul}(A) ) )
\to
V^{\bul}(A)\,,
\end{equation}
where $V^{\bul}(A)$ is equipped with the
vanishing differential.

Thus it suffices to connect the DG Gerstenhaber
algebras
$$
\Om_{\Ger} (B_{\Ger^{\vee}}( V^{\bul}(A) ) )
$$
and
$$
\Om_{\Ger} (\bbF_{\Ger^{\vee}}(\Cbu(A)), M)
$$
by a chain of quasi-isomorphisms.
For this purpose we will use the
sub $\Ger^{\vee}$-coalgebra
$$
\Xi(A)= \bbF_{\La^{2}\cocomm}(
\bbF_{\La\colie^+} (A, \Si\Der(A)) )
$$
of $\bbF_{\Ger^{\vee}}(\Cbu(A))$
introduced in equation (\ref{DT})\,.

Notice that, $\Xi(A)$ is also a sub
$\Ger^{\vee}$-coalgebra
of $\bbF_{\Ger^{\vee}}(V^{\bul}(A))$\,.
On the other hand,
due to remark $1$ after theorem \ref{adjun1}
$$
B_{\Ger^{\vee}}(V^{\bul}(A)) =
\bbF_{\Ger^{\vee}}(V^{\bul}(A))
$$
as $\Ger^{\vee}$-coalgebras.
Thus the natural question arises of whether
$\Xi(A)\subset \bbF_{\Ger^{\vee}}(V^{\bul}(A))$
 is stable under the
codifferential
\begin{equation}
\label{d-VA}
d_{V^{\bul}(A)} : \bbF_{\Ger^{\vee}}(V^{\bul}(A))
\to \bbF_{\Ger^{\vee}}(V^{\bul}(A))
\end{equation}
of $B_{\Ger^{\vee}}(V^{\bul}(A))$.

The following proposition gives the
positive answer to this question:
\begin{pred}
\label{sub-DG}
Let $d_{V^{\bul}(A)}$ be the codifferential
of the DG $\Ger^{\vee}$-coalgebra
$B_{\Ger^{\vee}}(V^{\bul}(A))$\,. Then
$$
d_{V^{\bul}(A)}\Big|_{\Xi(A)} \subset \Xi(A)\,.
$$
\end{pred}
{\bf Proof.}
According to the construction
(see p. $26$ in \cite{GJ}) of the functor
$B_{\cO^{\vee}}$ of the adjunction (\ref{adjun})
the proposition follows easily from the obvious
inclusions
$$
A \cdot \Der(A) \subset \Der(A)\,,
$$
$$
[\Der(A), A]_{SN} \subset A\,, \qquad
[\Der(A), \Der(A)]_{SN} \subset \Der(A)\,,
$$
where $\cdot$ and $[\,,\,]_{SN}$ denote the
product and the Lie bracket in $V^{\bul}(A)$,
respectively. $\Box$

We denote the embedding of $\Xi(A)$
(\ref{DT}) into
the DG $\Ger^{\vee}$-coalgebra
$B_{\Ger^{\vee}}(V^{\bul}(A))$ by $\io$
\begin{equation}
\label{io}
\io : \Xi(A) \hookrightarrow B_{\Ger^{\vee}}(V^{\bul}(A))\,.
\end{equation}
We will also use the same notation
$d_{V^{\bul}(A)}$ for the restriction of
the codifferential (\ref{d-VA}) to $\Xi(A)$\,.

Similarly to proposition \ref{sub-DG},
theorem \ref{ona} shows that
the natural embedding of $\Ger^{\vee}$-coalgebras
\begin{equation}
\label{sigma}
\si : \Xi(A) \hookrightarrow
\bbF_{\Ger^{\vee}}(\Cbu(A))
\end{equation}
is also compatible with the codifferentials
$d_{V^{\bul}(A)}$ and $M$\,, where
the codifferential $M$ on
$\bbF_{\Ger^{\vee}}(\Cbu(A))$ is the one
corresponding to the homotopy Gerstenhaber
algebra structure on $\Cbu(A)$\,.

Thus we get the following pair of morphisms
of DG $\Ger^{\vee}$-coalgebras:
\begin{equation}
\label{para}
(\bbF_{\Ger^{\vee}}(\Cbu(A)), M)
\,\stackrel{\si}{\leftarrow}\,
(\Xi(A), d_{V^{\bul}(A)})
\stackrel{\io}{\to}
B_{\Ger^{\vee}}( V^{\bul}(A) )
\end{equation}
Applying to (\ref{para}) the functor $\Om_{\Ger}$
(\ref{adjun}) we get the pair of morphisms of
DG Gerstenhaber algebras:
\begin{equation}
\label{para1}
\Om_{\Ger}( \bbF_{\Ger^{\vee}}(\Cbu(A)), M)
\,\stackrel{\Om_{\Ger} (\si) }{\leftarrow}\,
\Om_{\Ger}(\Xi(A), d_{V^{\bul}(A)})
\stackrel{ \Om_{\Ger} (\io) }{\to}
\Om_{\Ger} (B_{\Ger^{\vee}}( V^{\bul}(A) ) )\,.
\end{equation}
Combining this pair
with (\ref{chi-C}) and (\ref{eta-V}) we get
the following diagram of (homotopy) morphisms
of (homotopy) Gerstenhaber
algebras
\begin{equation}
\label{sigma-iota1}
\begin{array}{ccc}
\Cbu(A) & \stackrel{\chi^{\Cbu(A)}}{\brarrow} &
\Om_{\Ger} (\,\bbF_{\Ger^{\vee}}(\Cbu(A)), M \, ) \\[0.3cm]
 ~ & ~ &  \uparrow^{\Om_{\Ger}(\si)}  \\[0.3cm]
 ~ & ~ & \Om_{\Ger} ( \Xi(A), d_{V^{\bul}(A)} )    \\[0.3cm]
 ~ & ~ &  \downarrow^{\Om_{\Ger}(\io)}  \\[0.3cm]
V^{\bul}(A) & \,\stackrel{\eta_{\Ger}}{\leftarrow} \,
 & \Om_{\Ger}(\, B_{\Ger^{\vee}}(V^{\bul}(A)) \, ) \,,
\end{array}
\end{equation}
where $d_{V^{\bul}(A)}$ is the restriction of
the codifferential (\ref{d-VA}) to $\Xi(A)$\,.

Notice that, in diagram (\ref{sigma-iota1})
$\Om_{\Ger}(\si)$, $\Om_{\Ger}(\io)$, and
$\eta_{\Ger}$ (\ref{counit1}) are honest morphisms of
DG Gerstenhaber algebras and $\chi^{\Cbu(A)}$
is a homotopy morphism.

The main result of this paper can be
formulated as follows:
\begin{teo}
\label{domik}
Let $A$ be a regular commutative (associative) algebra
(with unit) over a field $\bbK$ of characteristic
zero. If $\Cbu(A)$ is the normalized
Hochschild cochain complex of $A$,
\begin{equation}
\label{M1}
M : \bbF_{\Ger^{\vee}}(C^{\bul}(A))
\to  \bbF_{\Ger^{\vee}}(C^{\bul}(A))
\end{equation}
is the homotopy Gerstenhaber algebra structure
on $C^{\bul}(A)$ induced by the action of the
little disk operad, $\Xi(A)$ is the
$\Ger^{\vee}$-coalgebra introduced in
(\ref{DT}), $\Om$ and
$B$ are the adjoint functors (\ref{adjun}),
and $\io$, $\si$ are embeddings
defined in (\ref{io}) and
(\ref{sigma})
then (\ref{sigma-iota1})
is a chain of quasi-isomorphisms
of (homotopy) Gerstenhaber algebras.
\end{teo}
{\bf Proof.}
The map $\eta_{\Ger}$ in (\ref{sigma-iota1})
is a quasi-isomorphism by theorem \ref{adjun1} and
$\chi^{\Cbu(A)}$ in (\ref{sigma-iota1}) is
a quasi-isomorphism by proposition \ref{pro-chi}\,.
To prove that $\Om_{\Ger}(\io)$ is
a quasi-isomorphism we show that so is
the composition $\eta_{\Ger} \circ \Om_{\Ger}(\io)$\,.
Then the fact that $\Om_{\Ger}(\si)$
is a quasi-isomorphism will follow from the
Hochschild-Kostant-Rosenberg theorem for
$\Cbu(A)$ \cite{HKR}\,.

Let us consider the composition
\begin{equation}
\label{nu}
\nu = \eta_{\Ger} \circ \Om_{\Ger}(\io)
:  \Om_{\Ger} ( \Xi(A), d_{V^{\bul}(A)} )  \to
V^{\bul}(A)\,.
\end{equation}
This is a map of DG Gerstenhaber algebras where
the algebra $V^{\bul}(A)$ (\ref{polyvec}) carries the
vanishing differential.

Since $\Der(A)$ is a module over $A$ we may
regard the pair $(A, \Si\Der(A))$ as a $\comm^+$-algebra.
Therefore, due to
equation \ref{comm-Kos}, the functor
$B_{\La\colie^+}$ of the adjunction (\ref{adjun})
gives us the pair
\begin{equation}
\label{BBB}
(B_{\La\colie}(A), B_{\La\colie^+}(A, \Si\Der(A))^+)
\end{equation}
of the DG $\La\colie$
coalgebra $B_{\La\colie}(A)$ and
its DG comodule $B_{\La\colie^+}(A, \Si\Der(A))^+$\,.
Notice that, as a  $\La\colie$  comodule
$B_{\La\colie^+}(A, \Si\Der(A))^+)$ is a cofree
comodule over $B_{\La\colie}(A)$ cogenerated by
$\Der(A)$ placed in degree $1$\,.

It is convenient to denote by
$B_{\La\colie^+}(A, \Si\Der(A))$ the direct
sum
\begin{equation}
\label{oplus}
B_{\La\colie^+}(A, \Si\Der(A)) =
 B_{\La\colie}(A)\oplus B_{\La\colie^+}(A, \Si\Der(A))^+\,.
\end{equation}
It is easy to see that (\ref{oplus}) is naturally
a DG $\La\colie$-coalgebra.
Furthermore, since the direct sum $A \oplus \Der(A)$
carries a Lie bracket which is compatible
with the product on $A$ and the $A$-module
structure on $\Der(A)$ the graded vector space
$B_{\La\colie^+}(A, \Si\Der(A))$ is 
equipped\footnote{See subsection \ref{razdel}
in which we give a detailed explanation of
how the functor $B_{\La\colie}$ can be extended to 
a functor from the category of Gerstenhaber 
algebras to the category of DG Lie bialgebras 
with Lie brackets of degree $-1$ and Lie 
cobrackets of degree $1$.}
with $\La\lie$-algebra structure which is
compatible with the
DG $\La\colie$-coalgebra.
Due to this fact we can apply
to $B_{\La\colie^+}(A, \Si\Der(A))$ the
functor $B_{\La^2\cocomm}$ and get
the following DG $\Ger^{\vee}$-coalgebra
$$
B_{\La^2\cocomm}(B_{\La\colie^+}(A, \Si\Der(A)))
$$

It is not hard to see that
\begin{equation}
\label{Xi}
(\Xi(A), d_{V^{\bul}(A)}) \cong
 B_{\La^2\cocomm}(B_{\La\colie^+}(A, \Si\Der(A)))
\end{equation}
as DG $\Ger^{\vee}$ coalgebras, where
$d_{V^{\bul}(A)}$ is, as above, the restriction of
the codifferential of $B_{\Ger^{\vee}}(V^{\bul}(A))$
to $\Xi(A)$\,.

Thus our purpose is to prove that
the map $\nu$ (\ref{nu})
\begin{equation}
\label{nuu}
\nu :  \Om_{\comm}\, \Om_{\La\lie}\,
 B_{\La^2\cocomm}\, B_{\La\colie^+}(A, \Si\Der(A))
\to
V^{\bul}(A)\,.
\end{equation}
is a quasi-isomorphism of DG Gerstenhaber algebras.

To do this we notice that applying the
functor $\Om_{\comm}$ to the DG 
$\La\colie$-coalgebra  
$B_{\La\colie^+}(A, \Si\Der(A))$ (\ref{oplus})
we get the following DG Gerstenhaber algebra
$$
\Om_{\comm}(B_{\La\colie^+}(A, \Si\Der(A)))
$$
with the obvious commutative product and the 
Lie bracket induced from 
$B_{\La\colie^+}(A, \Si\Der(A))$\,. 

Furthermore, there is canonical
map (of DG Gerstenhaber algebras)
\begin{equation}
\label{nu1}
\nu_1 : \Om_{\comm} \, \Om_{\La\lie}\,
B_{\La^2\cocomm} \, B_{\La\colie^+} (A, \Si\Der(A))
\to
\Om_{\comm} \, B_{\La\colie^+} (A, \Si\Der(A))\,,
\end{equation}
which is obtained by applying the functor
$\Om_{\comm}$ to the counit (\ref{counit1})
$$
\eta_{\La\lie} :  \Om_{\La\lie}\,
B_{\La^2\cocomm} \, B_{\La\colie^+} (A, \Si\Der(A))
\to B_{\La\colie^+} (A, \Si\Der(A))\,.
$$

On the other hand we can define the following map
of (DG) Gerstenhaber algebras
\begin{equation}
\label{nu2}
\nu_2 : \Om_{\comm} \, B_{\La\colie^+} (A, \Der(A))
\to V^{\bul}(A)\,.
\end{equation}
As a map of commutative algebras $\nu_2$ is
defined by its restriction to the space \\
$B_{\La\colie^+} (A, \Der(A))$
of generators of  $\Om_{\comm} \, B_{\La\colie^+} (A, \Der(A))$\,.
Namely, we define $\nu_2$ by setting
\begin{equation}
\label{nu2def}
\nu_2 \Big|_{B_{\La\colie^+} (A, \Der(A))}
= i \circ u_{A\, \oplus\, \Si \Der(A)}\,,
\end{equation}
where  $u_{A\, \oplus\, \Si \Der(A)}$ is the
corestriction (\ref{corest}) from
$B_{\La\colie^+} (A, \Der(A))$ to the space of its
cogenerators $A \oplus \Si \Der(A)$ and
$i$ is the natural
embedding\footnote{Recall that
$V^{\bul}(A)= S_A(\Si\Der(A))$ (\ref{polyvec}).}
$$
i :  A \oplus \Si \Der(A) \hookrightarrow
V^{\bul}(A)\,.
$$
It is not hard to check that the map $\nu_2$
(\ref{nu2}) defined by (\ref{nu2def}) is
compatible both with the differentials and the
Lie brackets.

Next, one can show that
the map $\nu$ (\ref{nuu}) is the
composition
\begin{equation}
\label{compos}
\nu= \nu_2 \circ \nu_1
\end{equation}
of the maps (\ref{nu1}) and (\ref{nu2})\,.

In order to prove that $\nu_1$ (\ref{nu1}) is
a quasi-isomorphism we introduce the following 
increasing filtration on the 
space $B_{\La\colie^+} (A, \Si\Der(A))$
\begin{equation}
\label{filtr}
\begin{array}{c}
\displaystyle
F^{0} B_{\La\colie^+} (A, \Si\Der(A)) 
\subset
\dots
\subset
F^{m} B_{\La\colie^+} (A, \Si\Der(A)) \subset \\[0.3cm]
\displaystyle
F^{m + 1}  B_{\La\colie^+} (A, \Si\Der(A)) \subset
\dots \subset B_{\La\colie^+} (A, \Si\Der(A))\,,
\end{array}
\end{equation}
where the subspace 
$F^m B_{\La\colie^+} (A, \Si\Der(A))$
is spanned by vectors
$$
(\la \otimes p)^{S_n}\,, \qquad n \le m + 1\,,
$$
where $\la\in \La\colie(n)$ and 
$p$ is a degree $n$ monomial in
$\bigotimes(A\oplus \Der(A))$ in which
an element of $\Der(A)$ can appear at most once.
This filtration is obviously  
compatible with the algebraic structures 
on $B_{\La\colie^+} (A, \Si\Der(A))$. Futhermore, 
both the codifferential and the   
cobracket on $B_{\La\colie^+} (A, \Si\Der(A))$ 
lower the filtration by $1$\,. 
Hence, on the associated graded space 
$\Gr (B_{\La\colie^+} (A, \Si\Der(A))) $ 
only the structure of a 
$\La\lie$-algebra survives, while the 
cobracket and the codifferential vanish.

The filtration (\ref{filtr}) canonically extends to 
exhausting increasing filtrations
on the complexes
$$
\Om_{\comm} \, \Om_{\La\lie}\,
B_{\La^2\cocomm} \, B_{\La\colie^+} (A, \Si\Der(A))
$$
and
$$
\Om_{\comm} \, B_{\La\colie^+} (A, \Si\Der(A))\,.
$$
The map $\nu_1$ (\ref{nu1}) is obviously compatible with
the resulting filtrations. Moreover, since the 
cobracket in  $B_{\La\colie^+} (A, \Si\Der(A))$ 
lowers the filtration by $1$ the
associated graded complexes are
$$
\bbF_{\comm} \, \Om_{\La\lie}\,
B_{\La^2\cocomm} \, \Gr ( B_{\La\colie^+} (A, \Si\Der(A)) )
$$
and
$$
\bbF_{\comm} \, \Gr( B_{\La\colie^+} (A, \Si\Der(A)) ) \,,
$$
respectively, and the associated graded
morphism 
$$
\Gr(\nu_1) = \bbF_{\comm} (\eta_{\La\lie})\,,
$$
where 
$$
\eta_{\La\lie} : \Om_{\La\lie}\,
B_{\La^2\cocomm} 
\, \Gr ( B_{\La\colie^+} (A, \Si\Der(A)) )
\to  \Gr (B_{\La\colie^+} (A, \Si\Der(A)) )\,,
$$
is the counit (\ref{counit1})
of the adjunction (\ref{adjun}), and 
$\Gr ( B_{\La\colie^+} (A, \Si\Der(A)))$ is viewed 
as the $\La\lie$-algebra.

By theorem \ref{adjun1} $\eta_{\La\lie}$
is a quasi-isomorphism. On the other hand,
the functor $\bbF_{\comm}$ is exact since
the underlying field $\bbK$ has characteristic
zero. Thus the map $\nu_1$ (\ref{nu1}) is
indeed a quasi-isomorphism.

Let us prove that the map $\nu_2$ (\ref{nu2}) is
a quasi-isomorphism.

Due to theorem \ref{adjun1} and proposition
\ref{comm-Koszul} we have the
quasi-isomorphism
\begin{equation}
\label{prrr}
\eta_{\comm}:
\Om_{\comm}(B_{\La\colie} (A)) \to A
\end{equation}
of DG commutative algebras (with
$A$ carrying the vanishing differential)
and the quasi-isomorphism
\begin{equation}
\label{prrr1}
\eta_{\comm^+}:
\Om_{\comm}(B_{\La\colie} (A))\otimes
B_{\La\colie^+} (A, \Si\Der(A))^+  \to
\Si\Der(A)
\end{equation}
of the corresponding (DG) modules
over the (DG) commutative algebras
$\Om_{\comm}(B_{\La\colie} (A))$ and
$A$\,.

Since $A$ is regular $\Der(A)$ is a flat module
over $A$. Hence the maps (\ref{prrr}) and
(\ref{prrr1})
induce the following quasi-isomorphism of
DG commutative algebras
\begin{equation}
\label{q-iso1}
\hat{\eta} :
S_{\Om_{\comm}(B_{\La\colie} (A))}(\Om_{\comm}(B_{\La\colie} (A))\otimes
B_{\La\colie^+} (A, \Si\Der(A))^+ ) \to
S_A(\Si\Der(A))
\end{equation}
On the other hand, it is obvious that
the DG commutative algebra
$$
S_{\Om_{\comm}(B_{\La\colie} (A))}(\Om_{\comm}(B_{\La\colie} (A))\otimes
B_{\La\colie^+} (A, \Si\Der(A))^+ )
$$
is isomorphic to
$$
\Om_{\comm} (B_{\La\colie^+}(A, \Si\Der(A)) )\,,
$$
where the DG $\La\colie$-coalgebra
$B_{\La\colie^+}(A, \Si\Der(A))$ is defined
in equation (\ref{oplus}).

Thus $\hat{\eta}$ (\ref{q-iso1}) is in fact
a quasi-isomorphism
\begin{equation}
\label{q-iso11}
\hat{\eta} :
\Om_{\comm} (B_{\La\colie^+}(A, \Si\Der(A)))
 \to
S_A(\Si\Der(A))\,.
\end{equation}
It is not hard to see that $\hat{\eta}$ coincides
with the map $\nu_2$ (\ref{nu2}). Hence,
$\nu_2$ is a quasi-isomorphism.

Due to equation (\ref{compos}) the map
$\nu$ (\ref{nu}) is a quasi-isomorphism.
Hence, so is $\Om_{\Ger}(\io)$\,.

It remains to prove that $\Om_{\Ger}(\si)$ in
(\ref{sigma-iota1}) is a quasi-isomorphism.
To prove this fact we notice that due
to the regularity of $A$ the cohomology of
$\Cbu(A)$ is generated by classes of $A\subset
\Cbu(A) $ and $\Der(A)\subset \Cbu(A)$\,.
The restriction of $\Om_{\Ger}(\si)$ to these
representatives in $\Om_{\Ger} (\Xi(A), d_{V^{\bul}(A)}) $
and $\Om_{\Ger} (\bbF_{\Ger^{\vee}}(\Cbu(A)),  M) $
gives the identity map. Hence, $\Om_{\Ger}(\si)$ is a
quasi-isomorphism and the theorem is proved. $\Box$

\section{Generalizations and applications}
One can notice that the proof of theorem
\ref{domik} is based on the flatness of
the module $\Der(A)$ over $A$ and the fact that
the cohomology of $\Cbu(A)$ is generated by
classes of $A\subset \Cbu(A)$ and
$\Der(A) \subset \Cbu(A)$\,.

This observation allows us to reformulate the
proof of theorem \ref{domik} and get the following
result:
\begin{teo}
\label{smooth}
If $\cM$ is a smooth real manifold
and $D^{\bul}(\cM)$ is the
graded vector space of polydifferential
operators on $\cM$ then the homotopy Gerstenhaber
algebra on $D^{\bul}(\cM)$ induced by
the map (\ref{vesch}) is formal. $\Box$
\end{teo}

Let $X$ be a smooth complex manifold
(resp. a smooth algebraic variety over a
field $\bbK$ of characteristic zero) and
$D^{\bul}_X$ be the sheaf of holomorphic
(resp. algebraic) polydifferential operators
on $X$. Due to the map (\ref{vesch}) $D^{\bul}_X$
is a sheaf of homotopy Gerstenhaber algebras.
Using the argument similar to the argument \cite{Vey}
for smooth real manifold it is not hard to show
that the sheaf of cohomology of $D^{\bul}_X$ is
the exterior algebra $\wedge^{\bul} T_X$ of
the tangent sheaf $T_X$. The product in
$\wedge^{\bul} T_X$ is the ordinary exterior product
and the bracket is the Schouten-Nijenhuis bracket
(see eq. ($3.20$) page $50$ in \cite{thesis}).

Replacing $A$ by $\cO_X$, $\Der(A)$ by
$T_X$, and $\Cbu(A)$ by $D^{\bul}_X$ in
diagram (\ref{sigma-iota1}) we get a chain
of homotopy morphisms between sheaves
of (homotopy) Gerstenhaber algebras connecting
the sheaves $D^{\bul}_X$ and $\wedge^{\bul} T_X$\,.
Generalizing the proof of theorem \ref{domik}
to this chain of morphisms we get the
following theorem:
\begin{teo}
\label{complex-alg}
The sheaf $D^{\bul}_X$ of
homotopy Gerstenhaber algebras is homotopy
equivalent to the sheaf of its
cohomology $\wedge^{\bul} T_X$\,. $\Box$
\end{teo}

Let us apply this theorem
and prove the last claim
in paper \cite{K}:
\begin{cor}[M. Kontsevich, \cite{K}, claim 8.4]
\label{ona1}
Let $X$ be either a smooth complex manifold
or a smooth algebraic variety (over a field
$\bbK$ of characteristic zero). Let
$\wedge^{\bul} T_X$ be the exterior algebra
of the tangent sheaf $T_X$ and
$\D$ be the diagonal in $X\times X$\,.
Then the graded commutative algebras
$$
H^{\bul}(X,  \wedge^{\bul} T_X)
$$
and
$$
Ext^{\bul}_{X\times X}(\cO_{\D}, \cO_{\D})
$$
are isomorphic.
\end{cor}
{\bf Proof.}
First, due to theorem $0.3$ in \cite{Ye1}
the ring
$Ext^{\bul}_{X\times X}(\cO_{\D}, \cO_{\D})$
is isomorphic to the hypercohomology
\begin{equation}
\label{HD-X}
\bbH^{\bul} (D^{\bul}_X)
\end{equation}
of the sheaf $D^{\bul}_X$ of
polydifferential operators on $X$\,.
The product on (\ref{HD-X}) is induced
by the ordinary cup-product of
polydifferential operators
\begin{equation}
\label{cup}
P_1 \cup P_2 (a_1, \dots, a_{k_1 + k_2}) =
P(a_1, \dots, a_{k_1})\, Q(a_{k_1+1}, \dots ,
a_{k_1 + k_2})\,,
\end{equation}
where $P_1$ (resp. $P_2$) is a polydifferential
operator of degree $k_1$ (resp. $k_2$).

The degree zero binary operation
of the homotopy
Gerstenhaber algebra on $D^{\bul}_X$
induces another commutative product on
the graded vector space
(\ref{HD-X}). Let us denote this product
by $\cup'$.

By theorem \ref{complex-alg} the rings
$$
(\bbH^{\bul} (D^{\bul}_X), \cup' )
$$
and
$$
H^{\bul}(X,  \wedge^{\bul} T_X)
$$
are isomorphic.
Therefore, it suffices to prove that
$\cup'$ coincides
with the product induced by (\ref{cup})
on (\ref{HD-X})\,.

It is clear from the construction \cite{M-Smith}
of J. E. McClure and J. H. Smith that the degree
zero binary operation $\cdot$ of the homotopy
Gerstenhaber algebra on $D^{\bul}_X$ is homotopic to
the cup-product (\ref{cup}). In other words, for
every pair of polydifferential operators
$P_1$ and $P_2$
\begin{equation}
\label{cup-dot}
P_1 \cdot P_2 = P_1 \cup P_2 + \pa  \Psi(P_1,P_2)
+ \Psi (\pa P_1, P_2) + (-1)^{|P_1|}
\Psi(P_1, \pa P_2)\,,
\end{equation}
where
$\pa$ is the Hochschild coboundary
operator, $|P_1|$ stands for the degree of
$P_1$ and $\Psi$ is a binary operation of
degree $-1$\,.

Let $(\check{C}^{\bul}(D^{\bul}_X), \cpa + \pa)$ be the
Cech complex of the sheaf $D^{\bul}_X$\,.
Let $P$ and $Q$ be cocycles in $\cCb(D^{\bul}_X)$
representing cohomology
classes $[P]$ and $[Q]$ in (\ref{HD-X}).

Equation (\ref{cup-dot}) implies that
$$
P \cdot Q = P \cup Q + \pa  \Psi(P,Q)
+ \Psi (\pa P, Q) + (-1)^{|P|}
\Psi(P, \pa Q)\,.
$$
On the other hand $P$ and $Q$ are cocycles
with respect to the differential $\cpa + \pa$\,.
Hence,
$$
P \cdot Q = P \cup Q + \pa  \Psi(P,Q)
- \Psi (\cpa P, Q) - (-1)^{|P|}
\Psi(P, \cpa Q)
$$
or equivalently
$$
P \cdot Q = P \cup Q + (\pa + \cpa) \Psi(P,Q)\,.
$$
Thus the product $\cup'$ coincides
with the product induced by (\ref{cup})
on (\ref{HD-X}) and the corollary follows. $\Box$

~\\
{\bf Remark.} Notice that (\ref{cup}) equips
the sheaf $D^{\bul}_X$ with a structure of a sheaf
of DG associative algebras. In particular, the
Cech complex $\check{C}^{\bul}(D^{\bul}_X)$ is
naturally a DGA. It is not hard to prove that
the graded algebra of cohomology
of $\check{C}^{\bul}(D^{\bul}_X)$ is
the graded (commutative) algebra
$H^{\bul}(X, \wedge^{\bul}\, T_X)$.
However, the question of whether
$\check{C}^{\bul}(D^{\bul}_X)$ is quasi-isomorphic
to its cohomology is much more subtle \cite{K1}.
As far as we know the answer to this question
is not yet found.

\section{Concluding remarks}
We would like to mention papers
\cite{CF} and \cite{Tomsk} in which
a super version of Kontsevich's formality
theorem was applied to quantum reduction.
In both of these papers the authors found
an obstruction (a quantum anomaly in the terminology
of \cite{Tomsk}) to quantization of
the quotient space. Notice that, the proof
of theorem \ref{domik}
admits an obvious generalization to the case of
a graded ring $A$\,. In particular, one can
easily extend a super version of Kontsevich's
formality quasi-isomorphism constructed in
\cite{CF} to a quasi-isomorphism of the
corresponding homotopy Gerstenhaber algebras.
It would be interesting to analyze the obstruction
of A.S. Cattaneo, G. Felder, S.L. Lyakhovich, and
A.A. Sharapov using the chain of quasi-isomorphisms
(\ref{sigma-iota1})\,.

We hope that the chain of quasi-isomorphisms
(\ref{sigma-iota1}) could be used to prove
a part of Caldararu's conjecture
formulated\footnote{See conjecture 5.2 in
\cite{Cald1}.} in \cite{Cald1}.

\end{document}